\documentclass[10pt]{article}
\usepackage{amsmath}
\usepackage{amssymb}
\usepackage{amsthm}
\usepackage{amscd}
\usepackage{amsfonts}
\usepackage[left=1in, top=1in,bottom=1in]{geometry}
\usepackage{sectsty}
\usepackage{xcolor}

\let\oldmarginpar\marginpar
\renewcommand\marginpar[1]{\-\oldmarginpar[\raggedleft\footnotesize #1]%
{\raggedright\footnotesize\sffamily #1}}

\allsectionsfont{\sffamily}

\usepackage{latexsym}
\usepackage{enumerate}
\usepackage{graphicx}
\usepackage{mathtools}
\newcommand{\bx}{\mathbf{x}}

\newcommand{\bTheta}{\boldsymbol{\Theta}}
\newcommand{\bF}{\mathbf{F}}
\newcommand{\bM}{\mathbf{M}}
\newcommand{\bX}{\mathbf{X}}
\newcommand{\bY}{\mathbf{Y}}
\newcommand{\bZ}{\mathbf{Z}}
\newcommand{\by}{\mathbf{y}}

\newcommand{\bmu}{\boldsymbol{\mu}}

\newcommand{\cB}{\mathcal{B}}

\newcommand{\cN}{\mathcal{N}}

\newcommand{\dB}{d\mathbf{B}}

\newcommand{\bbR}{\mathbb{R}}
\newcommand{\bbE}{\mathbb{E}}

\newcommand{\bbP}{\mathbb{P}}
\newcommand{\tr}{\!\top\!}
\newcommand{\scal}[2]{\left\langle{#1},{#2}\right\rangle}
\newcommand{\nor}[1]{\left\|{#1}\right\|}

\theoremstyle{plain}
\newtheorem{thm}{Theorem}

\theoremstyle{definition}

\newtheorem{rem}{Remark}

\DeclareMathOperator{\trace}{tr}

\title{\bfseries\sffamily
Wasserstein Contraction of Stochastic Nonlinear Systems
}

\author{
Jake Bouvrie%
\thanks{Laboratory for Computational and Statistical Learning, Massachusetts Institute of Technology, {\tt\small jvb@csail.mit.edu}}\\
\and
Jean-Jacques Slotine%
\thanks{Nonlinear Systems Laboratory, Massachusetts Institute of Technology,
        {\tt\small jjs@mit.edu}}\\
}
\date{\vspace{-5ex}}

\begin{document}
\maketitle
%\thispagestyle{empty}
%\pagestyle{empty}

%%%%%%%%%%%%%%%%%%%%%%%%%%%%%%%%%%%%%%%%%%%%%%%%%%%%%%%%%%%%%%%%%%%%%%%%%%%%%%%%
\begin{abstract}
We suggest that the tools of contraction analysis for deterministic systems can be applied towards
studying the convergence behavior of stochastic dynamical systems in the Wasserstein metric. In particular, we consider the case of Ito diffusions with identical dynamics but different initial condition distributions. If the drift term of the diffusion is contracting, then we show that the Wasserstein distance between the laws of any two solutions can be bounded by the Wasserstein distance between the initial condition distributions. Dependence on initial conditions exponentially decays in time at a rate governed by the contraction rate of the noise-free dynamics. The choice of the Wasserstein metric affords several advantages: it captures the underlying geometry of the space, can be efficiently estimated from samples, and advances a viewpoint which begins to bridge the gap between somewhat distinct areas of the literature.
\end{abstract}

%We suggest that the tools of contraction analysis for deterministic systems can be applied to the study %of coupled stochastic dynamical systems, with Wasserstein distance and the theory of optimal transport %serving as the key intermediary. If the drift term in an Ito diffusion is contracting, then we show %that the Wasserstein distance between the laws of two solutions to an Ito SDE exponentially forgets %initial conditions and can be bounded in the steady state by a constant factor governed by the %magnitude of the noise and the contraction rate of the noise-free dynamics. This viewpoint attempts to %contribute a step towards extending the contraction framework to stochastic systems, and bridging the %gap between somewhat distinct areas of the literature.

%%%%%%%%%%%%%%%%%%%%%%%%%%%%%%%%%%%%%%%%%%%%%%%%%%%%%%%%%%%%%%%%%%%%%%%%%%%%%%%%
\section{Introduction}
In this brief note we point out an application of contraction analysis for deterministic systems~\cite{Lohmiller98} to a class of It\^o stochastic differential equations, extending the mean-square stochastic contraction framework in~\cite{Pham09} towards considering Wasserstein distance~\cite{Villani03} as the metric for studying convergence. Contraction theory aims to provide convenient global exponential convergence results for coupled systems, and a Wasserstein viewpoint begins to illuminate the deeper connections between contraction and optimal transport. The Wasserstein metric is also natural for studying stochastic systems because it captures underlying geometry in an intuitive sense, via the earth mover's interpretation~\cite{Villani03}, while mean-square or total-variation distances for instance do not take geometry into account. Deterministic contracting systems also enjoy rich structural properties which make the analysis of complex networks of dynamical elements more tractable~\cite{SlotineModular03}, and although it is not our focus here, we expect that similar properties can be shown to hold for classes of stochastic systems. 

Stability and convergence of distributions governed by dynamical systems is in general well-studied, and the development that follows builds upon, and complements, the analyses of contracting stochastic systems discussed in~\cite{Pham09} and~\cite{Tabareau10,TabareauRDS}; work from the more optimal-transport or PDE oriented literature~\cite{Sturm05,Natile11,Bolley12}; and recent optimization perspectives~\cite{Raginsky2017,WibisonoColt2018} which consider gradient flows as a special case. The contraction-theoretic viewpoint we adopt here attempts to highlight useful perspectives from contraction in the context of stochastic systems, and bridge the gap between overlapping yet somewhat distinct areas of the literature.

\section{Background}
We briefly recall the basic definitions from contraction theory and optimal transport needed to support the development that follows below. We refer the reader to the references for a complete treatment and further background.

\subsection{Contraction Analysis}\label{sec:contraction-review}
Contraction analysis as originally proposed in~\cite{Lohmiller98,WangSlotine05} pertains to general,  deterministic systems of the form $\dot{\bx} = f(x,t)$, with $f:\bbR^{d}\times[0,\infty)\to\bbR^{d}$, and makes a broad statement about the conditions under which trajectories of such systems can be expected to globally converge exponentially fast to a single nominal trajectory. Here, $f$ is assumed to be real and smooth enough so that any required derivatives exist and are continuous. The basic result is as follows:

\begin{thm}[{\rm Contraction~\cite{Lohmiller98}}]\label{thm:contraction}
Consider the system $\dot{\bx} = f(\bx,t)$ in $\bbR^d$, with $f$ smooth and nonlinear. If there exists a uniformly positive definite metric
$\bM(\bx,t) = \bTheta(\bx,t)^{\tr}\bTheta(\bx,t)$
such that the associated generalized Jacobian
\[
\bF = \left(\dot{\bTheta} + \bTheta\frac{\partial f}{\partial \bx}\right)\bTheta^{-1}
\]
is uniformly negative definite, then all system trajectories converge exponentially to a single trajectory with convergence rate $|\sup_{\bx,t}\lambda_{\rm max}(\bF)|>0$ where $\lambda_{\rm max}$ is the largest eigenvalue of the symmetric part of $\bF(\bx,t)$, and the system is said to be contracting.
\end{thm}

If $f$ is a conservative vector field so that $f=\nabla\phi$ for some potential function $\phi:\bbR^d\to\bbR$, then $f$ contracting implies that $-\phi$ is strongly convex. If $\phi\in C^2$, then an equivalent second order condition for contraction is that $-\nabla^2\phi \succeq \beta I$ for some $\beta>0$.

\subsection{Wasserstein Distance}\label{sec:wasser-review}
The Wasserstein distance is a natural metric for comparing probability distributions that nicely captures the underlying geometry of the space it's defined on. It is often convenient to work with in applications, and intuitive to grasp via the optimal transport ``earth mover's'' interpretation. We borrow from the theory of optimal transport developed by Villani~\cite{Villani03}, where an extensive exposition can be found.

If $\mu,\nu$ are two probability measures on $\bbR^d$ with bounded second moments, the 2-Wasserstein distance between them is defined as:
\[
W_{2}(\mu,\nu) := \inf\bigl(\bbE\|X-Y\|^2\bigr)^{1/2}
\]
where $\|\cdot\|$ is the $\ell_2$ norm on $\bbR^d$ and the infimum is taken over the set of all joint measures on $\bbR^d \times \bbR^d$ with marginals $X\sim\mu$ and $Y\sim\nu$ (i.e. the set of couplings of $\mu,\nu$). Convergence in $W_2$ also implies convergence of the first two moments.
%The 2-Wasserstein distance metrizes the space of measures with bounded second moments, and convergence %in $W_2$ implies ordinary weak convergence of measures plus convergence of the first two moments.

\section{Stochastic Contraction in Wasserstein Distance}
We now show that a simple contraction result similar in spirit to Theorem~\ref{thm:contraction} holds in a Wasserstein sense for certain It\^o stochastic differential equations. However, now, the quantity for which we will  establish a suitable form of contraction is the distribution describing solutions to a stochastic differential equation (SDE) (or one can equivalently look at the flow of the Fokker-Planck equation counterpart to a given SDE). Our approach borrows from and is inspired by~\cite{Pham09}, where contraction for a class of stochastic systems is established in a mean-square sense. 
%We will consider for simplicity a uniform contraction metric, which without loss of generality, may be %chosen to be the identity matrix.

%We first define a preliminary Lemma.
%\begin{lem}\label{lem:cont-lem}
%Let $f$ be a continuously differentiable vector-valued function on $\bbR^d$. If the Jacobian of $f$ is strictly negative definite, i.e. if $J_f(\bx) \preceq -\beta I$ with $\beta >0$ uniformly for all $\bx$, then 
%$$\bigl(f(\bx) - f(\by)\bigr)^{\tr}(\bx-\by)\leq -\beta\nor{\bx-\by}^2.$$
%\end{lem}
%\begin{rem}
%If $f$ is a conservative vector field so that $f=\nabla\phi$ for some potential function $\phi:\bbR^d\to\bbR$, then $f$ contracting implies that $-\phi$ is strongly convex. If $\phi\in C^2$, then \mbox{$-\nabla^2\phi \succeq \beta I$} is an equivalent second order condition.
%\end{rem}
%\begin{proof}
%By the mean value theorem for vector valued functions, 
%$$\bigl(f(\bx) - f(\by)\bigr)^{\tr}(\bx-\by) = (\bx-\by)^{\tr}J_f(\bz)(\bx-\by)$$
%with $\bz=\eta\bx + (1-\eta)\by$ for some $\eta\in(0,1)$. The Lemma then follows using that $\bu^{\tr}J_f\bu\leq -\beta\nor{\bu}^2$ for any vector $\bu$.
%\end{proof}

Our main result says that if a sufficiently well-behaved It\^o stochastic differential equation has a contracting drift function, then the SDE can be said to be exponentially contracting in the 2-Wasserstein distance. Here we restrict ourselves to potentially time-varying contraction metrics $\bM(t)$; see~\cite{Pham13} for a discussion of state-dependent metrics\footnote{Note that the contraction metric in Theorem~\ref{thm:contraction} should not be confused with the Wasserstein metric; these are distinct quantities.}.
\begin{thm}[{\rm Stochastic Contraction in $W_2$}]\label{thm:stoch_cont}
Consider the It\^o stochastic differential equation
\begin{equation}\label{eqn:was-sde}
d\bX_t = f(\bX_t)dt + \sigma(\bX_t,t)\dB_t, \quad t\in[0,T], \bX_0\sim\mu_0
\end{equation}
with $\bbE\nor{\bX_0}^2<\infty$, $\mathbf{B}\in\bbR^d$ a standard Brownian motion, 
$\bX_t\in\bbR^d$, $f:\bbR^d\to\bbR^d$, $\sigma:\bbR^d\times[0,T]\to\bbR^{d\times d}$,
 and assume that the following conditions hold:
\begin{enumerate}[(i)]
\item The drift function $f$ is contracting in some uniformly positive definite metric $\bM(t)$
with contraction rate bounded by a finite constant $\beta>0$, where $\bM(t)$ satisfies 
$\bx^{\tr}\bM(t)\bx \geq \alpha\nor{\bx}^2$ with $\alpha>0$ for all $\bx, t$. (\cite[Condition H1']{Pham09})
\item For all $\bx, t\in[0,T]$,
\[
\trace\bigl(\sigma(\bx,t)^{\tr}\bM(t)\sigma(\bx,t)\bigr) \leq C_{\sigma}
\]\label{lab:cond_noisevar}
\item The diffusion matrix $\sigma(\bx,t)\sigma(\bx,t)^{\tr}$ satisfies 
\[
\scal{\sigma(\bx,t)\sigma(\bx,t)^{\tr}\by}{\by}\geq c\nor{\by}^2, \quad c >0
\] uniformly for all $\bx, \by\in\bbR^d$, $t\in[0,T]$. \label{lab:cond_ellip}
\item There exist constants $K_1,K_2 >0$ such that $\forall t\in[0,T], \bx,\by\in\bbR^d$,
\begin{align*}
\nor{f(\bx)} + \nor{\sigma(\bx,t)}_{F} &\leq K_1(1+\nor{\bx})\\
\nor{f(\bx)-f(\by)} + \nor{\sigma(\bx,t)-\sigma(\by,t)}_{F} &\leq K_2\nor{\bx-\by} .
\end{align*}\label{lab:cond_lip}
\end{enumerate}
Then the system~\eqref{eqn:was-sde} is {\em stochastically contracting} in the sense that, for any pair of solutions $\bX_t, \bY_t$ to~\eqref{eqn:was-sde} with respective laws $\mu_t, \nu_t$,
\begin{equation}\label{eqn:wasser_contr}
W_2(\mu_t, \nu_t)\leq 
\alpha^{-1/2}\bigl(e^{-\beta t}W_2(\mu_0, \nu_0) + \sqrt{C_{\sigma}/\beta}\bigr), \quad t\in[0,T].
\end{equation}
\end{thm}
\begin{rem}
By choosing $\nu_0$ equal to the stationary distribution of~\eqref{eqn:was-sde}, we see that the above result suggests exponentially fast convergence of $\mu_t$ to equilibrium in the space of probability measures. However, as we will discuss below, the effect of the noise does not vanish in~\eqref{eqn:wasser_contr} as solutions become close due to the assumption of distinct, independent noise processes driving $\bX_t$ and $\bY_t$ and an approach that passes through a mean-square calculation (see also~\cite{Pham09} Sec. 2.2.1, 2.3.2 for a related discussion).
%Note that while $W_2$ gives a weak convergence, it can be strengthened further when combined with smoothness bounds~\cite{Villani03}.
\end{rem}
\begin{proof}
%The first portion of the proof through~\eqref{eqn:gron_bound} essentially follows~\cite[Theorem 2]{Pham09}, which we include here for completeness. 
Condition~(\ref{lab:cond_lip}) is standard and guarantees existence and uniqueness of solutions to~\eqref{eqn:was-sde}, while condition~(\ref{lab:cond_ellip}) ensures existence of a unique invariant distribution. Suppose $\bX_t, \bY_t$ are two solutions to~\eqref{eqn:was-sde} driven by independent noise processes, starting from initial conditions independent of the noise and governed by the distributions $\mu_0$ and $\nu_0$ respectively:
\begin{align*}
d\bX_t &= f(\bX_t)dt + \sigma(\bX_t,t)\dB_t^{(1)}, &\quad \bX_0 &\sim\mu_0\\
d\bY_t &= f(\bY_t)dt + \sigma(\bY_t,t)\dB_t^{(2)}, &\quad \bY_0 &\sim\nu_0 \;.
\end{align*}
Let $\mu_t,\nu_t$ denote the laws of $\bX_t,\bY_t$ respectively.
With the definitions and conditions above, we can apply~\cite[Theorem 3]{Pham09} to obtain the mean-square relation
\begin{equation}\label{eqn:gron_bound}
\bbE\nor{\bX_t-\bY_t}^2 
%&\leq e^{-2\beta t}\bbE\nor{\bX(0)-\bY(0)}^2 + \frac{C_{\sigma}}{\beta}(1-e^{-2\beta t})\\
\leq \frac{1}{\alpha}\left(e^{-2\beta t}\bbE\nor{\bX_0-\bY_0}^2 + \frac{C_{\sigma}}{\beta}\right) .
\end{equation}

From~\eqref{eqn:gron_bound} we can show contraction in the 2-Wasserstein distance by essentially just taking the infimum over all couplings of $\mu_0$ and $\nu_0$. Define the shorthand $v(z):=\nor{x-y}^2$ for $z=(x,y)$, and let $\bZ_t:=(\bX_t,\bY_t)$. 
Let $P:[0,T]\times\bbR^{2d}\times\cB(\bbR^{2d})\to\bbR_{+}$ denote the transition function of the Markov process $\bZ_t$, and recall that $P(t,z,B)=\bbP(\bZ_t\in B | \bZ_0=z)$ almost surely. Now suppose $\pi_0^*$ is an optimal coupling of $\mu_0$ and $\nu_0$ so that $W_2^2(\mu_0, \nu_0)=\int v(z)\pi_0^*(dz)$. Finally, let $\bbE_{\pi_0^*}$ denote expectation with respect to the product measure formed from the measure $\pi_0^*$ on $\bZ_0$ and the (independent) Wiener measure induced by the noise on the SDE describing $\bZ_t$. Then,
\begin{align*}
\bbE_{\pi_0^*}[v(\bZ_t)] &= \int_{\bbR^{2d}}\int_{\bbR^{2d}} v(z)P(t,z_0,dz)\pi_0^*(dz_0) \\
&= \int_{\bbR^{2d}} \bbE[v(\bZ_t) ~|~ \bZ_0=z_0] \pi_0^*(dz_0) \\
&\leq \int_{\bbR^{2d}} \alpha^{-1}\Bigl(e^{-2\beta t}v(z_0) + C_{\sigma}/\beta\Bigr)\pi_0^*(dz_0) \\
&= \alpha^{-1}\bigl(e^{-2\beta t}W_2^2(\mu_0, \nu_0) + C_{\sigma}/\beta\bigr)
\end{align*}
where the third line follows from the second substituting in Equation~\eqref{eqn:gron_bound} to estimate the inner expectation. Defining the measure $\pi_t(B):=\int P(t, z_0, B)\pi_0^*(dz_0)$, we may view $\bbE_{\pi_0^*}[v(\bZ_t)]= \int v(z)\pi_t(dz)$ as an expectation of $v(\bZ_t)$ with respect to the coupling $\pi_t$. Since $\pi_t$ is not necessarily an optimal coupling,  $W^2_2(\mu_t, \nu_t) \leq \bbE_{\pi_0^*}[v(\bZ_t)]$. Taking square roots then gives the Theorem.
% SIMPLE ENDING: optimize the RHS over all couplings to get the tightest bound, which is with Wasser
% on the RHS.
%From the definition of the $2$-Wasserstein distance, $W^2_2(\mu_t, \nu_t) \leq \bbE\nor{\bX_t-\bY_t}^2$.
%Taking the infimum on both sides over all couplings between $\bX_0$ and $\bY_0$ with marginals $\mu_0$ %and $\nu_0$, it follows that
%\[
%W^2_2(\mu_t, \nu_t) \leq e^{-2\beta t}W^2_2(\mu_0, \nu_0) + \frac{C_{\sigma}}{\beta} .
%\]
%Taking square roots gives the Theorem.
\end{proof}

\subsection{Discussion}
The Theorem shows that if the drift of an It\^o diffusion is contracting, the distributions describing solutions to the SDE~\eqref{eqn:was-sde} given different initial conditions forget those initial conditions exponentially fast, and can be shown to be close in the Wasserstein metric.\footnote{It may be the case that exponential convergence of solutions in the Wasserstein metric is {\em equivalent} to contraction of the noise-free dynamics $f$; see for instance~\cite{Sturm05,Natile11}.} Closeness here is governed by the constant $\sqrt{C_{\sigma}/\alpha\beta}$ which depends solely on the noise amplitude and the contraction rate of the noise-free dynamics, however this term results from the assumption of independent noise processes driving different solutions to the SDE, and the use of the mean-square estimate~\eqref{thm:stoch_cont}. The assumption of independent noise sources is an important one, as we would like contraction to say something useful about e.g. coupled physical systems such as oscillators, each of which would be subject to its own noise, or to study the relationship between noisy and noise-free solutions of a system. It is possible, however, that a less conservative residual bound could be obtained by considering instead the associated Fokker-Planck equation for~\eqref{eqn:was-sde}.
% But it is possible that the noise variance term unavoidably introduced by the mean-square bound~\eqref{thm:stoch_cont} can be eliminated by considering an alternative development (for instance, one that examines the Fokker-Planck equation for~\eqref{eqn:was-sde}).

While the mean-square relationship~\eqref{eqn:gron_bound} is used in the course of establishing the Wasserstein bound, the two do not provide identical information, and Wasserstein distances have some advantages. The Wasserstein distance is stronger and more revealing in that it takes into account the underlying geometry of the space (while mean-square or total variation do not), giving a more intuitive notion of distance between distributions. Wasserstein distances can also be estimated efficiently from samples; we can approximate a Wasserstein distance between two distributions with the Wasserstein distance between their empirical distributions (see e.g.~\cite{Fournier2015}). In applications, the quantity $W_2(\mu_0, \nu_0)$ would benefit from both of these observations thereby providing a potentially more illuminating statement about the behavior of a stochastic system.

Finally, it is helpful to study the behavior of SDEs from an optimal transportation standpoint
in the context of optimization related applications. The Langevin dynamics in particular have been studied as a means to understand sampling-based algorithms, where sampling is viewed as optimization in the space of probability measures~\cite{Raginsky2017,WibisonoColt2018}. This raises natural connections between contraction and tools used within the optimization and machine learning communities, and  underscores the potential role of contraction in understanding the behavior of optimization algorithms more generally (see e.g.~\cite{Wensing18}).

\subsection{Example: Ornstein-Uhlenbeck System}
We will show Theorem~\ref{thm:stoch_cont} at work in a simplified but illustrative setting where the key quantities can be calculated analytically, and where we can appraise tightness of the bound given in the Theorem by comparing it to an exact Wasserstein calculation. Consider the linear SDE in $\bbR^d$
\begin{equation}\label{eqn:ou-sde}
d\bX_t = A(\bmu-\bX_t)dt + \sigma\dB_t
\end{equation}
where $A\in\bbR^{d\times d}$ is a strictly positive definite matrix and $\sigma > 0$ is scalar.  The solution to this system $\bX_t$ is an Ornstein-Uhlenbeck (OU) process. 

The noise-free dynamics $\dot{\bx} = A(\bmu-\bx)$ can be seen to be contracting in the identity metric ($\alpha=1$) with rate $\underline{\lambda}(A):=\lambda_{\rm min}(A+A^{\tr})/2 > 0$ given the assumption that $A$ is strictly positive definite, while we have  $C_{\sigma}=d\sigma^2$ in this simple setting. The other conditions $(ii)$-$(iv)$ of Theorem~\ref{thm:stoch_cont} can also be shown to hold given the assumption $A$ is strictly positive definite and $\sigma > 0$. The result of Theorem~\ref{thm:stoch_cont} then gives us a bound of 
$\sqrt{d\sigma^2/\underline{\lambda}(A)}$ on the 2-Wasserstein distance between the laws of two OU process solutions to~\eqref{eqn:ou-sde} after an exponential transient of rate $\underline{\lambda}(A)$. 

We can proceed with solving~\eqref{eqn:ou-sde} explicitly to confirm the claimed exponential convergence, and compare the steady-state bound to the known fact that $W_2$ actually goes to zero. Applying the variation of constants method,
\begin{equation}
\label{eqn:ou-soln}
\bX_t = e^{-At}\bX_0 + \bigl(I - e^{-At}\bigr)\bmu  +
\sigma\int_0^t e^{-A(t-s)}\dB_s
\end{equation}
and $\bX_t$ is seen to be a Gaussian process that can be characterized entirely by its time-dependent mean and covariance, $\bX_t\sim\cN\bigl(\bmu(t),\Sigma(t)\bigr)$. A straightforward calculation gives 
\begin{align}
\bmu(t) :&=\bbE[\bX_t] = e^{-At}\bbE[\bX_0] + \bigl(I - e^{-At}\bigr)\bmu\label{eqn:w-mean-process}\\
\Sigma(t) :&= \bbE\left[\bigl(\bX_t -\bbE\bX_t\bigr)\bigl(\bX_t -\bbE\bX_t\bigr)^{\!\tr}\right]\notag\\
 &= e^{-At}\Sigma(0)e^{-A^{\tr}t}
+ \sigma^2(A+A^{\tr})^{-1}\bigl(I - e^{-(A+A^{\tr})t}\bigr) .\notag
\end{align}
It is clear that the stationary distribution has the form  $\bX_{\infty}:=\lim_{t\to\infty}\bX_t\sim\cN\bigl(\bmu, (\sigma^2/4)(A^{-1}+A^{-\tr})\bigr)$, convergence to this steady-state is exponential, and  $W_2$ must go to 0 as $t\to\infty$ for any pair of laws respectively governing any pair of solutions with different initial condition distributions. A quick check that $W_2\to 0$ is possible here, using that the 2-Wasserstein distance between two Gaussian distributions $w_1=\cN\bigl(\mu_1,\Sigma_1\bigr)$ and $w_2=\cN\bigl(\mu_2,\Sigma_2\bigr)$ on $\bbR^d$ is given explicitly by
\[
W^2_2(w_1,w_2) = \|\mu_1 - \mu_2\|_2^2 + \trace\bigl(\Sigma_1 + \Sigma_2 
		- 2(\Sigma_2^{1/2}\Sigma_1\Sigma_2^{1/2})^{1/2}\bigr) \;.
\]
So clearly $W_2\to 0$ as $t\to\infty$ for two Gaussian process solutions to~\eqref{eqn:ou-sde} of the form~\eqref{eqn:w-mean-process}, as desired.

%\subsection{Coupled Noisy Oscillators}
%\subsection{Networked Stochastic Nonlinear Systems}

%%%%%%%%%%%%%%%%%%%%%%%%%%%%%%%%%%%%%%%%%%%%%%%%%%%%%%%%%%%%%%%%%%%%%%%%%%%%%%%%
%\section{ACKNOWLEDGMENTS}
%The authors gratefully acknowledge.

%%%%%%%%%%%%%%%%%%%%%%%%%%%%%%%%%%%%%%%%%%%%%%%%%%%%%%%%%%%%%%%%%%%%%%%%%%%%%%%%

\bibliographystyle{plain}
\bibliography{wasser_contr}

\end{document}